\documentclass[12pt]{amsart}
\usepackage{amsmath,amssymb,amsfonts}
\begin{document}
\theoremstyle{plain}
\newtheorem{thm}{Theorem}[subsection]
\newtheorem{lem}[thm]{Lemma}
\newtheorem{cor}[thm]{Corollary}
\newtheorem{prop}[thm]{Proposition}
\newtheorem{rem}[thm]{Remark}
\newtheorem{defn}[thm]{Definition}
\newtheorem{ex}[thm]{Example}

\numberwithin{equation}{subsection}
\newcommand{\mc}{\mathcal}
\newcommand{\mb}{\mathbb}
\newcommand{\surj}{\twoheadrightarrow}
\newcommand{\inj}{\hookrightarrow}
\newcommand{\zar}{{\rm zar}}
\newcommand{\an}{{\rm an}} 
\newcommand{\red}{{\rm red}}
\newcommand{\codim}{{\rm codim}}
\newcommand{\rank}{{\rm rank}}
\newcommand{\Ker}{{\rm Ker \ }}
\newcommand{\Pic}{{\rm Pic}}
\newcommand{\Div}{{\rm Div}}
\newcommand{\Hom}{{\rm Hom}}
\newcommand{\im}{{\rm im}}
\newcommand{\Spec}{{\rm Spec \,}}
\newcommand{\Sing}{{\rm Sing}}
\newcommand{\Char}{{\rm char}}
\newcommand{\Tr}{{\rm Tr}}
\newcommand{\Gal}{{\rm Gal}}
\newcommand{\Min}{{\rm Min \ }}
\newcommand{\Max}{{\rm Max \ }}
\newcommand{\sA}{{\mathcal A}}
\newcommand{\sB}{{\mathcal B}}
\newcommand{\sC}{{\mathcal C}}
\newcommand{\sD}{{\mathcal D}}
\newcommand{\sE}{{\mathcal E}}
\newcommand{\sF}{{\mathcal F}}
\newcommand{\sG}{{\mathcal G}}
\newcommand{\sH}{{\mathcal H}}
\newcommand{\sI}{{\mathcal I}}
\newcommand{\sJ}{{\mathcal J}}
\newcommand{\sK}{{\mathcal K}}
\newcommand{\sL}{{\mathcal L}}
\newcommand{\sM}{{\mathcal M}}
\newcommand{\sN}{{\mathcal N}}
\newcommand{\sO}{{\mathcal O}}
\newcommand{\sP}{{\mathcal P}}
\newcommand{\sQ}{{\mathcal Q}}
\newcommand{\sR}{{\mathcal R}}
\newcommand{\sS}{{\mathcal S}}
\newcommand{\sT}{{\mathcal T}}
\newcommand{\sU}{{\mathcal U}}
\newcommand{\sV}{{\mathcal V}}
\newcommand{\sW}{{\mathcal W}}
\newcommand{\sX}{{\mathcal X}}
\newcommand{\sY}{{\mathcal Y}}
\newcommand{\sZ}{{\mathcal Z}}
\newcommand{\A}{{\mathbb A}}
\newcommand{\B}{{\mathbb B}}
\newcommand{\C}{{\mathbb C}}
\newcommand{\D}{{\mathbb D}}
\newcommand{\E}{{\mathbb E}}
\newcommand{\F}{{\mathbb F}}
\newcommand{\G}{{\mathbb G}}
\renewcommand{\H}{{\mathbb H}}
\newcommand{\I}{{\mathbb I}}
\newcommand{\J}{{\mathbb J}}
\newcommand{\M}{{\mathbb M}}
\newcommand{\N}{{\mathbb N}}
\renewcommand{\P}{{\mathbb P}}
\newcommand{\Q}{{\mathbb Q}}
\newcommand{\R}{{\mathbb R}}
\newcommand{\T}{{\mathbb T}}
\newcommand{\U}{{\mathbb U}}
\newcommand{\V}{{\mathbb V}}
\newcommand{\W}{{\mathbb W}}
\newcommand{\X}{{\mathbb X}}
\newcommand{\Y}{{\mathbb Y}}
\newcommand{\Z}{{\mathbb Z}}


\catcode`\@=11
%
%
\def\opn#1#2{\def#1{\mathop{\kern0pt\fam0#2}\nolimits}} 
\def\bold#1{{\bf #1}}%
\def\underrightarrow{\mathpalette\underrightarrow@}
\def\underrightarrow@#1#2{\vtop{\ialign{$##$\cr
 \hfil#1#2\hfil\cr\noalign{\nointerlineskip}%
 #1{-}\mkern-6mu\cleaders\hbox{$#1\mkern-2mu{-}\mkern-2mu$}\hfill
 \mkern-6mu{\to}\cr}}}
\let\underarrow\underrightarrow
\def\underleftarrow{\mathpalette\underleftarrow@}
\def\underleftarrow@#1#2{\vtop{\ialign{$##$\cr
 \hfil#1#2\hfil\cr\noalign{\nointerlineskip}#1{\leftarrow}\mkern-6mu
 \cleaders\hbox{$#1\mkern-2mu{-}\mkern-2mu$}\hfill
 \mkern-6mu{-}\cr}}}
%
%

%
\def\:{\colon}
\let\oldtilde=\tilde
\def\tilde#1{\mathchoice{\widetilde{#1}}{\widetilde{#1}}%
{\indextil{#1}}{\oldtilde{#1}}}
\def\indextil#1{\lower2pt\hbox{$\textstyle{\oldtilde{\raise2pt%
\hbox{$\scriptstyle{#1}$}}}$}}
\def\pnt{{\raise1.1pt\hbox{$\textstyle.$}}}
%

%
\let\amp@rs@nd@\relax
\newdimen\ex@
\ex@.2326ex
\newdimen\bigaw@
\newdimen\minaw@
\minaw@16.08739\ex@
\newdimen\minCDaw@
\minCDaw@2.5pc
\newif\ifCD@
\def\minCDarrowwidth#1{\minCDaw@#1}
\newenvironment{CD}{\@CD}{\@endCD}
\def\@CD{\def\A##1A##2A{\llap{$\vcenter{\hbox
 {$\scriptstyle##1$}}$}\Big\uparrow\rlap{$\vcenter{\hbox{%
$\scriptstyle##2$}}$}&&}%
\def\V##1V##2V{\llap{$\vcenter{\hbox
 {$\scriptstyle##1$}}$}\Big\downarrow\rlap{$\vcenter{\hbox{%
$\scriptstyle##2$}}$}&&}%
\def\={&\hskip.5em\mathrel
 {\vbox{\hrule width\minCDaw@\vskip3\ex@\hrule width
 \minCDaw@}}\hskip.5em&}%
\def\verteq{\Big\Vert&&}%
\def\noarr{&&}%
\def\vspace##1{\noalign{\vskip##1\relax}}\relax\let\amp@rs@nd@&\iffalse}\fi
 \CD@true\vcenter\bgroup\relax\let\\=\cr\iffalse}\fi\tabskip\z@skip\baselineskip20\ex@
 \lineskip3\ex@\lineskiplimit3\ex@\halign\bgroup
 &\hfill$\m@th##$\hfill\cr}
\def\@endCD{\cr\egroup\egroup}
%
\def\>#1>#2>{\amp@rs@nd@\setbox\z@\hbox{$\scriptstyle
 \;{#1}\;\;$}\setbox\@ne\hbox{$\scriptstyle\;{#2}\;\;$}\setbox\tw@
 \hbox{$#2$}\ifCD@
 \global\bigaw@\minCDaw@\else\global\bigaw@\minaw@\fi
 \ifdim\wd\z@>\bigaw@\global\bigaw@\wd\z@\fi
 \ifdim\wd\@ne>\bigaw@\global\bigaw@\wd\@ne\fi
 \ifCD@\hskip.5em\fi
 \ifdim\wd\tw@>\z@
 \mathrel{\mathop{\hbox to\bigaw@{\rightarrowfill}}\limits^{#1}_{#2}}\else
 \mathrel{\mathop{\hbox to\bigaw@{\rightarrowfill}}\limits^{#1}}\fi
 \ifCD@\hskip.5em\fi\amp@rs@nd@}
\def\<#1<#2<{\amp@rs@nd@\setbox\z@\hbox{$\scriptstyle
 \;\;{#1}\;$}\setbox\@ne\hbox{$\scriptstyle\;\;{#2}\;$}\setbox\tw@
 \hbox{$#2$}\ifCD@
 \global\bigaw@\minCDaw@\else\global\bigaw@\minaw@\fi
 \ifdim\wd\z@>\bigaw@\global\bigaw@\wd\z@\fi
 \ifdim\wd\@ne>\bigaw@\global\bigaw@\wd\@ne\fi
 \ifCD@\hskip.5em\fi
 \ifdim\wd\tw@>\z@
 \mathrel{\mathop{\hbox to\bigaw@{\leftarrowfill}}\limits^{#1}_{#2}}\else
 \mathrel{\mathop{\hbox to\bigaw@{\leftarrowfill}}\limits^{#1}}\fi
 \ifCD@\hskip.5em\fi\amp@rs@nd@}
%
%
\newenvironment{CDS}{\@CDS}{\@endCDS}
\def\@CDS{\def\A##1A##2A{\llap{$\vcenter{\hbox
 {$\scriptstyle##1$}}$}\Big\uparrow\rlap{$\vcenter{\hbox{%
$\scriptstyle##2$}}$}&}%
\def\V##1V##2V{\llap{$\vcenter{\hbox
 {$\scriptstyle##1$}}$}\Big\downarrow\rlap{$\vcenter{\hbox{%
$\scriptstyle##2$}}$}&}%
\def\={&\hskip.5em\mathrel
 {\vbox{\hrule width\minCDaw@\vskip3\ex@\hrule width
 \minCDaw@}}\hskip.5em&}
\def\verteq{\Big\Vert&}
\def\novarr{&}
\def\noharr{&&}
\def\SE##1E##2E{\slantedarrow(0,18)(4,-3){##1}{##2}&}
\def\SW##1W##2W{\slantedarrow(24,18)(-4,-3){##1}{##2}&}
\def\NE##1E##2E{\slantedarrow(0,0)(4,3){##1}{##2}&}
\def\NW##1W##2W{\slantedarrow(24,0)(-4,3){##1}{##2}&}
\def\slantedarrow(##1)(##2)##3##4{%
\thinlines\unitlength1pt\lower 6.5pt\hbox{\begin{picture}(24,18)%
\put(##1){\vector(##2){24}}%
\put(0,8){$\scriptstyle##3$}%
\put(20,8){$\scriptstyle##4$}%
\end{picture}}}
\def\vspace##1{\noalign{\vskip##1\relax}}\relax\let\amp@rs@nd@&\iffalse}\fi
 \CD@true\vcenter\bgroup\relax\let\\=\cr\iffalse}\fi\tabskip\z@skip\baselineskip20\ex@
 \lineskip3\ex@\lineskiplimit3\ex@\halign\bgroup
 &\hfill$\m@th##$\hfill\cr}
\def\@endCDS{\cr\egroup\egroup}
%
\newdimen\TriCDarrw@
\newif\ifTriV@
\newenvironment{TriCDV}{\@TriCDV}{\@endTriCD}
\newenvironment{TriCDA}{\@TriCDA}{\@endTriCD}
\def\@TriCDV{\TriV@true\def\TriCDpos@{6}\@TriCD}
\def\@TriCDA{\TriV@false\def\TriCDpos@{10}\@TriCD}
\def\@TriCD#1#2#3#4#5#6{%
\setbox0\hbox{$\ifTriV@#6\else#1\fi$}
\TriCDarrw@=\wd0 \advance\TriCDarrw@ 24pt
\advance\TriCDarrw@ -1em
\def\SE##1E##2E{\slantedarrow(0,18)(2,-3){##1}{##2}&}
\def\SW##1W##2W{\slantedarrow(12,18)(-2,-3){##1}{##2}&}
\def\NE##1E##2E{\slantedarrow(0,0)(2,3){##1}{##2}&}
\def\NW##1W##2W{\slantedarrow(12,0)(-2,3){##1}{##2}&}
\def\slantedarrow(##1)(##2)##3##4{\thinlines\unitlength1pt
\lower 6.5pt\hbox{\begin{picture}(12,18)%
\put(##1){\vector(##2){12}}%
\put(-4,\TriCDpos@){$\scriptstyle##3$}%
\put(12,\TriCDpos@){$\scriptstyle##4$}%
\end{picture}}}
\def\={\mathrel {\vbox{\hrule
   width\TriCDarrw@\vskip3\ex@\hrule width
   \TriCDarrw@}}}
\def\>##1>>{\setbox\z@\hbox{$\scriptstyle
 \;{##1}\;\;$}\global\bigaw@\TriCDarrw@
 \ifdim\wd\z@>\bigaw@\global\bigaw@\wd\z@\fi
 \hskip.5em
 \mathrel{\mathop{\hbox to \TriCDarrw@
{\rightarrowfill}}\limits^{##1}}
 \hskip.5em}
\def\<##1<<{\setbox\z@\hbox{$\scriptstyle
 \;{##1}\;\;$}\global\bigaw@\TriCDarrw@
 \ifdim\wd\z@>\bigaw@\global\bigaw@\wd\z@\fi
 \mathrel{\mathop{\hbox to\bigaw@{\leftarrowfill}}\limits^{##1}}
 }
 \CD@true\vcenter\bgroup\relax\let\\=\cr\iffalse}\fi
 \tabskip\z@skip\baselineskip20\ex@
 \lineskip3\ex@\lineskiplimit3\ex@
 \ifTriV@
 \halign\bgroup
 &\hfill$\m@th##$\hfill\cr
#1&\multispan3\hfill$#2$\hfill&#3\\
&#4&#5\\
&&#6\cr\egroup%
\else
 \halign\bgroup
 &\hfill$\m@th##$\hfill\cr
&&#1\\%
&#2&#3\\
#4&\multispan3\hfill$#5$\hfill&#6\cr\egroup
\fi}
\def\@endTriCD{\egroup}


\title{ Algebraic Differential Characters }  
\author{ H\'el\`ene Esnault } 
\address{ Universit\"at Essen, FB6 Mathematik, 45 117 Essen, Germany}
\email{ esnault@uni-essen.de}
\thanks{ This work has been partly supported by the DFG Forschergruppe
''Arithmetik und Geometrie''}

\maketitle
\setcounter{section}{-1} 
\section{Introduction}
In \cite{CS}, Cheeger and Simons defined on a $\sC^{\infty}$
manifold $X$ a group of differential characters $\hat{H}^{2n}
(X, \R /\Z)$, which is an extension of the global $\R$ valued
closed forms of degree $2n$ having $\Z$ periods, by the group
$H^{2n-1} (X, \R/\Z)$. (In fact, they write $\hat{H}^{2n-1} (X,
\R/\Z)$ but the notation $2n$ rather than $(2n-1)$ fits better
with weights in algebraic geometry). Similarly, there is a group
of (complex) differential characters $\hat{H}^{2n} (X, \C/\Z)$,
presented as an extension of the global $\C$ valued closed forms
of degree $2n$ having $\Z$ periods, by the group $H^{2n-1} (X,
\C/\Z)$. The group $\hat{H}^{2n} (X, \R/\Z)$ (resp. $\hat{H}^{2n} (X,
\C/\Z)$) is also presented as an extension of the Betti
cohomology group $H^{2n} (X, \Z)$ by global $\R$ valued (resp.
$\C$ valued) differential forms of degree $2n-1$, modulo the
closed ones with $\Z$ periods.

They define a ring structure on $\hat{H}^{2 \bullet} (X, \R/\Z)$
(resp. $\hat{H}^{2 \bullet} (X, \C/\Z)$) and show the existence
of functorial and additive classes $\hat{c}_n (E, \nabla)_{\R}
\in \hat{H}^{2n} (X, \R/\Z)$ (resp. $c_n (E, \nabla) \in
\hat{H}^{2n} (X, \C/ \Z)$) for a $\sC^{\infty}$ bundle with a
connection $\nabla$, lifting the closed form $P_n (\nabla^2 ,
\ldots , \nabla^2)$ with $\Z$ periods, where $P_n$ is the
homogeneous symmetric invariant polynomial of degree $n$ which
is the $n$-th symmetric function in the entries on diagonal
matrices. When the connection $\nabla$ is flat, that is when
$\nabla^2 = 0$, then $\hat{c}_n (E, \nabla)_{\R} \in H^{2n-1}
(X, \R/\Z)$ (resp. $\hat{c}_n (E, \nabla) \in H^{2n-1} (X,
\C/\Z)$).

The aim of this article is to develop a similar construction
when $X$ is an algebraic smooth variety over a field $k$ and
$(E, \nabla)$ is an algebraic bundle with an algebraic
connection. We define a group of {\it algebraic differential
characters} $AD^n (X)$, which is an extension of global closed
algebraic forms of degree $2n$, whose cohomology class in
$\H^{2n} (X, \Omega^{\geq 2n}_{X})$ is algebraic, by the group
$\H^n (X, \sK_n \> d \log >> \Omega^{n}_{X} \to ...)$,
introduced and studied in \cite{EII}, \cite{BE}. When $k = \C$,
$AD^n (X)$ maps to the group of {\it analytic differential
characters} $D^n (X)$, defined as an extension of global
analytic forms of degree $2n$, with $\Z(n)$ periods, by
$H^{2n-1} (X_{\an}, \C /\Z (n))$. The group $D^n (X)$ maps to
the Deligne cohomology group $H^{2n}_{\sD} (X, \Z (n))$. (The
above description is valid when $X$ is proper, otherwise one has
to modify it a little bit). The group $AD^n (X)$ (resp. $D^n
(X)$) is also presented as an extension of the subgroup $\Ker
CH^n (X) \to \H^n (X, \Omega^{n}_{X} \to \ldots \to
\Omega^{2n-1}_{X})$ of the Chow group $CH^n (X)$ (resp. of the
subgroup $\Ker H^{2n} (X_{\an} , \Z (n)) \to H^n (X_{\an},
\Omega^{n}_{X} \to \ldots \to \Omega^{2n-1}_{X} ))$ of the Betti
cohomology group $H^{2n} (X_{\an} , \Z (n))$) by $\H^{n-1} (X,
\Omega^{n}_{X} \to \ldots \to \Omega^{2n-1}_{X})/H^{n-1}(X,
\sK_n)$ (resp.
$\H^{n-1} (X_{\an} , \Omega^{n}_{X} \to \ldots \to
\Omega^{2n-1}_{X} )/H^{2n-1}(X_{\an}, \Z(n))$). There is a ring structure on
$AD^{\bullet} (X), D^{\bullet} (X)$. If $(E, \nabla)$ is an
algebraic bundle with an algebraic connection, we define
functorial and additive classes $c_n (E, \nabla) \in AD^n (X)$,
lifting $P_n (\nabla^2, \ldots, \nabla^2)$ and the algebraic
Chern classes $c_n (E) \in CH^n (X)$. The group $AD^n (X)$ maps
to $H^0 (X, \Omega^{2n-1}_{X}/ d \Omega^{2n-2}_{X})$. The class
$c_n (E, \nabla)$ lifts the algebraic Chern-Simons class $w_n
(E, \nabla)$, related to the algebraic equivalence relation on
cycles, defined in \cite{BE}. When $\nabla$ is flat, then $c_n
(E, \nabla) \in \H^n (X, \sK_n \> d \log >> \Omega^{n}_{X} \to
\ldots )$ and is the class defined in \cite{EII}.

We give two constructions of the classes. The first one is a
{\it generalized splitting principle}. In \cite{EI} and
\cite{EII}, we had defined a modified splitting principle when
$\nabla$ is flat. However, our construction had the disavantage
to use the flag bundle of $E$, rather than simply its projective
bundle $\P(E)$, and to introduce a ``$\tau$ cohomology'' on the
flag bundle which is not a free module over the corresponding
``$\tau$ cohomology'' of the base $X$. We correct those two
points here by introducing a slightly more complicated ``$\tau$
cohomology'' on $\P (E)$, also defined when $\nabla$ is not
flat, which is free over the corresponding ``$\tau$ cohomology''
on $X$. The product structure is then naturally defined
following the recipe explained in \cite{B} and \cite{EV}. The
whole construction is now very closed to the construction of
Hirzebruch-Grothendieck for bundles without supplementary
structure.

The second construction relies on the generalized {\it Weil
algebra} defined by Beilinson and Kazhdan in \cite{BK}. To a
bundle $E$, they associate functorially a Weil algebra complex
$\Omega^{\bullet}_{X, E}$ together with a group $H^{2n} (X,
U_E (n))$, which is an extension of the homogeneous
symmetric invariant polynomials $P$ of degree $n$ 
mapping to 
$$
\text{Im} (H^{2n} (X_{\an}, \Z (n)) \to H^{2n} (X_{\an} , \C))$$
via the Weil homomorphism,
by $H^{2n-1} (X_{\an} , \C/\Z (n))$. By evaluating their
construction on the universal simplicial bundle, they defined
functorial classes 
$$c^{BK}_{n} (E) \in H^{2n} (X, U_E
(n)).$$ A connection $\nabla$ on $E$ defines a map $\nabla : H^n
(X, U_E (n)) \to D^n (X)$, and thereby classes in $D^n
(X)$.

We modify their construction to make it more algebraic. We
obtain in this way $c^{ABK}_{n} (E, \nabla) \in AD^n (X)$. We
show that $c_n (E, \nabla) = c^{ABK}_{n} (E, \nabla)$,
reflecting the fact that the splitting principle and the
universal bundle construction define the same Chern classes for
bundles without supplementary structure. 

\subsection{Acknowledgements:}

In the process of writing \cite{BE}, Spencer Bloch got convinced
of the existence of classes $c_n (E, \nabla) \in AD^n (X)$
lifting our algebraic Chern-Simons classes $w_n (E, \nabla)$,
and restricting to $c_n (E, \nabla)$ defined in \cite{EII} when
$\nabla^2 =0$. He explained to me the construction of
Beilinson-Kazhdan \cite{BK} mentioned above, out of which I
could perform the construction of $c^{ABK}_{n} (E, \nabla)$. It
is my pleasure to thank him for generously sharing his ideas
with me. I thank Eckart Viehweg for encouraging me at different
stages of this work, in particular for reminding me all the
tricks of the sign conventions in products we understood while
writing \cite{EV}. \section{Notations} \label{sec:not} 
\begin{enumerate}
\item[1)] $X$ is a smooth algebraic variety over a field $k$,
$D$ is a normal crossing divisor, $j: U = X - D \to X$ is the
embedding. 
\item[2)] $\lambda : X \to \bar{X}$ is a good compactification,
such that $\Delta = \bar{X} - X$ and $\bar{X} - U = \bar{D}$ are
normal crossing divisors, $\bar{j} = \lambda \circ j : U \to
\bar{X}$ is the embedding. 
\item[3)] $(\Omega^{\bullet}_{X} (\log D), \Omega^{\geq n}_{X}
(\log D))$ is the de Rham complex with logarithmic poles along
$D$, filtered by its stupid filtration, and $\Omega^{<n}_{X}
(\log D)$ is its quotient. 
\item[4)] $\tau_0 : \Omega^{\bullet}_{X} (\log D) \to
N^{\bullet}$ is a map of complexes, where $N^{\bullet}$ is a
differential graded algebra, such that the sheaves $N^n$ are
locally free, as well as $B^b = \Ker \Omega^{b}_{X} (\log D) \to
N^b$, such that $N^0 = \sO_X$, and that if $a$ is the smallest
degree for which $B^a \neq 0$, then $B^b = B^a \wedge
\Omega^{b-a}_{X} (\log D)$ for $b \geq a$ (see \cite{BE} (3.11)
and \cite{EI} (2.1)). 
\item[5)] $\sK_n$ is the image of the Zariski sheaf
$\sK^{M}_{n}$ of Milnor $K$ theory in $K^{M}_{n} (k (X))$, with
its $d \log $ map $\sK_n \> d \log >> \Omega^{\geq n}_{X} [n]$,
and its induced map $\tau_0 \circ d \log : \sK_n \to N^{\geq n}
[n]$. $N^{\infty} \sK_n$ is the complex $(\sK_n \to N^n \to
N^{n+1} \to ...)$. 
\item[6)] $(E, \nabla)$ is an algebraic bundle of rank $r$ with
an algebraic $\tau_0$ connection, also called a $N^1$ valued
connection 
$$
\nabla : E \to N^1 \otimes_{\sO_X} E, 
$$
that is a $k$ linear map verifying the Leibniz rule 
$$
\nabla (\lambda e) = \tau_0 \circ d (\lambda) \otimes e +
\lambda \nabla e.
$$
$\nabla^2$ is the $\sO_X$ linear map $\nabla \circ \nabla : E
\to N^2 \otimes E$, call the curvature of $\nabla$, with 
$$
\nabla : N^n \otimes E \to N^{n+1} \otimes E
$$
defined by the sign convention 
$$
\nabla ( \alpha \otimes e) = \tau_0 \circ d (\alpha ) \otimes e
+ (-1)^n \alpha \wedge \nabla e.
$$
If $\nabla^2 =0$, one defines the $\tau_0$ de Rham complex 
$$
(N^{\bullet} \otimes E, \nabla).
$$
\item[7)] If $k = \C$, $a: X_{\an} \to X_{\zar}$ is the identity
from $X$ endowed with the analytic topology to $X$ endowed with
the Zariski topology. For a sheaf $\sF$ on $X_{\zar}$, we denote
by $\sF_{\an}$ the sheaf $a^* \sF$. When the context is clear,
we still write $\sF$ for $\sF_{\an}$. 
\item[8)] $H^{P}_{\sD} (X, \Z (q))$ is the Deligne-Beilinson
cohomology of $X$. 
\item[9)] $\sF^a \sH^b$ is the Zariski sheaf associated to $F^a
H^{b}_{DR} (U)$ (the Hodge filtration on the de Rham
cohomology), $\sH^n (\C/\Z (m))$ to $H^n (U, \C/\Z (m))$,
$\sF^{a}_{\Z (c)} \sH^b = \Ker (\sF^a \sH^b \to \sH^b (\C/\Z
(c)))$, $\sH^{a}_{\sD} (b)$ to $H^{a}_{\sD} (U, \Z (b))$. 
\item[10)] For $G = GL (r)$, of $\sG = M (r)$, we denote by $P_n
\in S^n (\sG^*)^G$ the symmetric invariant polynomial of degree
$n$ which is the $n$-th symmetric function of the diagonal
entries on diagonal matrices. We denote by $G_e \subset G$
the matrices fixing the subspace spanned by the $r''$ first canonical
basis vectors, and by $\sG_e \subset \sG$ the corresponding algebra.  
\item[11)] If $E$ is a vector bundle, we denote by $\sG_E$ the
endomorphisms of $E$, and by $\sG^{*}_{E}$ its dual (of course
isomorph to $\sG_E$). If 
$$ (e) : 0 \to E'' \to E \to E' \to 0$$
is an exact sequence of bundles with 
$r'' = \text{rank} E'', r' = \text{rank} E'$, we denote by 
$\sG_{(E,e}) \subset \sG_E$ the endomorphisms respecting $(e)$
and by $\sG^*_{(E,e)}$ its dual.
\end{enumerate}


\section{Algebraic Differential Characters} \label{sec:adc} 

\begin{defn} \label{dfn:dc}
Let $\tau_0$ be as in \ref{sec:not}, 4). We denote by $\iota$
the natural embedding $N^{\geq 2n} [n] \to N^{\geq n} [n]$. We
define on $X$ the complex 
$$
C(n)_{\tau_0} = {\rm cone} \ (\sK_n \oplus N^{\geq 2n} [n] \> d
\log \oplus - \iota >> N^{\geq n} [n])[-1].
$$
If $k = \C$, we define  on $ \bar{X}_{\an}$
\begin{multline*}
DR (n)_{\tau_0} = {\rm cone} \ (\Omega^{<n}_{\bar{X}} (\log
\bar{D}) \to R \lambda_* (N^{\geq n} [1] ) [-1] \notag \\
= \sO_X \to \cdots \to \Omega^{n-1}_{\bar{X}} (\log \bar{D}) \to
\lambda_* N^n \to \lambda_* N^{n+1} \cdots 
\end{multline*}
$$
C(n)^{\an}_{\tau_0} = {\rm cone} \ (R \lambda_* \Z (n) \oplus R
\lambda_* (N^{\geq 2n}) [n] \> \epsilon \oplus -\iota >> DR
(n)_{\tau_0}) [-1]
$$
where $\epsilon$ is induced by the map 
$$
\epsilon: R \lambda_* \Z (n) \to \Omega^{\bullet}_{\bar{X}}
(\log \bar{D}). 
$$
We define 
$$
\begin{array}{ll} 
AD^{n}_{\tau_0} (X) = \H^n (X, C (n)_{\tau_0}) & \\
D^{n}_{\tau_0} (X) = \H^{2n} (X, C (n)^{\an}_{\tau_0}) &
\mbox{if} \ \ k = \C .
\end{array} 
$$
When $\tau_0 =$ identity and $D = \phi$, we simply write $AD^n
(X)$ and $D^n (X)$. 
\end{defn}
\begin{prop} \label{prop:diag} 
\ \\
\begin{enumerate}
\item[1)] The group $D^{n}_{\tau_0} (X)$ does not depend on
$\lambda $ (this justifies the notation). 
\item[2)] There are exact sequences
\begin{enumerate}
\item[i)] \begin{multline*} 
0 \to \H^n (X, N^{\infty} \sK_n) \to AD^{n}_{\tau_0} (X) \notag \\
\to \Ker (H^0 (X , N^{2n}_{cl}) \to \frac{\H^{2n} (X, N^{\geq
n})}{CH^n (X)} ) \to 0  \notag
\end{multline*}
\item[ii)] \begin{multline*}
0 \to \H^{2n} (X_{\an}, {\rm cone} \ (R \lambda_* \Z (n) \to DR
(n)_{\tau_0} ) [-1]) \\
\to D^{n}_{\tau_0} (X) \to \Ker (H^0 (X_{\an}, N^{2n}_{cl}) \\
\to \frac{\H^{2n} (\bar{X}_{\an} , DR (n)_{\tau_0})}{H^{2n}
(X_{\an}, \Z (n))} ) \to 0 
\end{multline*}
\end{enumerate}
For $\tau_0 =$ identity, $D = \phi$ and $X$ proper, this reads 
\begin{multline*}
0 \to H^{2n-1} (X_{\an} , \C / \Z (n)) \to D^n (X) \\
\to \Ker (H^0 (X_{\an} , \Omega^{2n}_{cl}) \to H^{2n} (X_{\an} ,
\C/\Z (n))) \to 0
\end{multline*} 
For $\tau_0 =$ identity, $D = \phi$ and $X$ non-proper,
then 
one has a splitting $\H^m(X_{\an}, DR(n)) = H^m_{DR}(X) \oplus R(m,n)$
and 
$$
\rm{Ker}\  D^n (X) \> \rm{restriction} >> 
\H^{2n} (X_{\an} , C (n)^{\an} |_{X_{\an}}) = R(2n-1,n) .
$$ 
\item[3)] If $k = \C$, there is a commutative diagram 
$$
\begin{CD}
AD^{n}_{\tau_0} (X) \> \alpha >> CH^n (X) \\
\V \psi_{\tau_0} VV \V V \psi V \\
D^{n}_{\tau_0} (X) \> \beta >> H^{2n}_{\sD} (X, \Z (n))
\end{CD}
$$
\end{enumerate}
\end{prop}

\begin{proof}
\ \\
\begin{enumerate}
\item[1)] As usual, if $\lambda$ and $\lambda'$ are two good
compactifications, one constructs a third one $\lambda_1$
dominating $\lambda$ and $\lambda'$, with $\sigma: \bar{X}_1 \to
\bar{X}$. One just has to compare the $\lambda$ and the
$\lambda_1$ constructions. Then $\sigma$ induces a map 
$$
\sigma^* : \H^m (\bar{X}, C (n)^{\an}_{\tau_0 , \lambda} ) \to
\H^m (\bar{X}_1 , C (n)^{\an}_{\tau_0, \lambda_1} )
$$
(where we put a $\lambda$ index to underline the dependance),
and from the cone definition of $C(n)^{\an}_{\tau_0 , \lambda}$,
one just has to see that 
$$\sigma^* : \H^m (\bar{X}, DR(n)_{\tau_0 , \lambda} ) \to
\H^m (\bar{X}_1, DR(n)_{\tau_0 , \lambda_1} ) 
$$
is an isomorphism. But one has a long exact sequence
\begin{multline*}
\to \H^m (X_{\an}, N^{\geq n}) \to \H^m (\bar{X}, DR (n)_{\tau_0
\lambda}) \\
\to H^{m}_{DR} (X) / F^n H^{m}_{DR} (X) \to ...
\end{multline*}
and $\sigma^*$ induces an isomorphism on the two terms $\H^m (X_{\an} ,
N^{\geq n})$ and $H^{m}_{DR} (X)/F^n H^{m}_{DR} (X)$. 
\item[2)] 
\begin{enumerate}
\item[i)] We just regard the exact triangle
$$
\to N^{\infty} \sK_n \to C (n)_{\tau_0} \to N^{\geq 2n} [n] \> [1] >>
$$
\item[ii)] Similarly we regard the exact triangle
\begin{multline*}
\to {\rm cone} \ (R \lambda_* \Z (n) \to DR (n)_{\tau_0} ) [-1] \\
\to C (n)^{\an}_{\tau_0} \to R \lambda_* N^{\geq 2n}_{\an} [n] \> [1] >>.
\end{multline*}
When $\tau_0 =$ identity and $D = \phi$, the maps
$$
\to \Omega^{\bullet}_{\bar{X}_{\an}} (\log \bar{D}) \to DR (n) \to R 
\lambda_* (\Omega^{\bullet}_{X_{\an}}) 
$$
define a splitting 
$$
\H^m (\bar{X}_{\an}, DR (n) ) = H^{m}_{DR} (X) \oplus R (m ,n),
$$
with $H^m(X_{\an}, \Z(n))$ mapping to $H^m_{DR}(X)$.
This shows that the restriction map is injective on the right
hand side
of the exact sequence ii) and that
the Kernel of the left hand side is $R(2n-1, n)$.
\end{enumerate}
\item[3)] $\alpha$ is induced by 
$$
C(n)_{\tau_0} \to \sK_n
$$
and $\beta$ by
$$
C(n)^{\an}_{\tau_0} \to {\rm cone} \ (R \lambda_* \Z (n) \to
\Omega^{<n}_{\bar{X}} (\log \bar{D})) [-1].
$$
The cycle map is compatible with restriction to open sets. From
the commutative diagrams of exact sequences 
$$
\begin{CD}
\noarr \noarr \H^{2n-1} (X, N^{\geq n} / N^{\geq 2n}) \=
\H^{2n-1} (X, N^{\geq n} / N^{\geq 2n} ) \noarr \\
\noarr \noarr \V VV \V VV \noarr \\
0 \>>> \Ker \>>> H \>>> AD^{n}_{\tau_0} (X) \>>> 0 \\
\noarr \verteq \V VV \V VV \noarr \\
0 \>>> \Ker \>>> CH^n (\bar{X}) \>>> CH^n (X) \>>> 0 \\
\noarr \noarr \V VV \V VV \noarr \\
\noarr \noarr \H^n (X, N^{\geq n} / N^{\geq 2n}) \= \H^{n} (X,
N^{\geq n} / N^{\geq 2n} ) \noarr 
\end{CD}
$$
\ \\ \\
$$
\begin{CD}
\noarr \noarr \H^{2n-1} (X_{\an}, N^{\geq n} / N^{\geq 2n}) \=
\H^{2n-1} (X_{\an} , N^{\geq n} / N^{\geq 2n} ) \\
\noarr \noarr \V VV \V VV \\
0 \>>> \Ker^{\an} \>>> H^{\an} \>>> D^{n}_{\tau_0} (X) \\
\noarr \verteq \V VV \V VV \\
0 \>>> \Ker^{\an} \>>> H^{2n}_{\sD} (\bar{X}, \Z (n)) \>>>
H^{2n}_{\sD} (X, \Z (n)) \\
\noarr \noarr \V VV \V VV \\
\noarr \noarr \H^{2n} (X_{\an}, N^{\geq n} / N^{\geq 2n}) \= \H^{2n} (X_{\an},
N^{\geq n} / N^{\geq 2n} )
\end{CD} 
$$
with 
\begin{gather}
H = \H^{n} (\bar{X}, \bar{C} (n)_{\tau_0} ) \notag \\
H^{\an} = \H^{2n} (\bar{X}, \bar{C} (n)^{\an}_{\tau_0} ) \notag \\
\bar{C} (n)_{\tau_0} = {\rm cone} \ (\sK_n \oplus R \lambda_*
(N^{\geq 2n} )[n] \>  d\log \oplus -\iota >> R \lambda_* (N^{\geq n}))
[-1] \notag \\
\bar{C} (n)^{\an}_{\tau_0} = {\rm cone} \ (\Z (n) \oplus R
\lambda_* (N^{\geq 2n} ) [n] \> \epsilon \oplus -\iota >> \bar{DR}
(n)_{\tau_0} ) [-1] \notag \\
\bar{DR} (n)_{\tau_0} = {\rm cone} \ (\Omega^{<n}_{\bar{X}} \to
R \lambda_* N^{\geq n} [1] ) [-1]. \notag 
\end{gather}
One sees that it is enough to lift $\psi_{\bar{X}} : CH^n
(\bar{X}) \to H^{2n}_{\sD} (\bar{X} , \Z (n))$ to
$\psi_{\bar{X}, \tau_0} : H \to H^{\an}$. One has 
$$
R^i a_* \bar{C} (n)^{\an}_{\tau_0} = \sH^{i-1} (\C / \Z (n)) \ \
\mbox{for} \ \ i \leq n -1.
$$
By the Bloch-Ogus vanishing theorem (\cite{BO}, (6.2)), this
implies 
\begin{multline*}
\H^{2n} (\bar{X} , \bar{C} (n)^{\an}_{\tau_0} ) = \\
\H^{2n}
(\bar{X}, {\rm cone}( \ 
\tau_{\leq n-1} R a_* \bar{C} (n)^{\an}_{\tau_0} \to
Ra_*\bar{C} (n)^{\an}_{\tau_0}) ) 
\end{multline*}
One has an exact sequence
$$
0 \to \sH^{n-1} (\C / \Z (n)) \to \sH^{n}_{\sD} (\Z (n)) \to
\sF^{n}_{\Z (n)} \sH^n \to 0 
$$
and the map of complexes 
\begin{multline*}
(\sH^{n}_{\sD} (\Z (n)) \to (a_* \lambda_* N)^{\geq n}/ (a_*
\lambda_* N)^{\geq 2n} [n]) \\
\to {\rm cone}( \ 
\tau_{\leq n-1} R a_* \bar{C} (n)^{\an}_{\tau_0} \to
Ra_*\bar{C} (n)^{\an}_{\tau_0})
\end{multline*}
factors through
\begin{multline*}
(\sF^{n}_{\Z (n)} \sH^n \to (a_* \lambda_* N)^{\geq n} / (a_*
\lambda_* N)^{\geq 2n} [n] ) \\
\> \sigma >> 
{\rm cone}( \ 
\tau_{\leq n-1} R a_* \bar{C} (n)^{\an}_{\tau_0} \to
Ra_*\bar{C} (n)^{\an}_{\tau_0}).
\end{multline*}
To obtain $\psi_{\bar{X}, \tau_0} : H \to H^{\an}$, we first map
$\sK_n \to (\lambda_* N)^{\geq n} / (\lambda_* N)^{\geq 2n}[n]$ to
$\sF^{n}_{\Z (n)} \sH^{n} \to (a_* \lambda_* N)^{\geq n} / (a_*
\lambda_* N)^{\geq 2n}[n]$ via the $d \log$ map on $\sK_n$ and the
change of topology map on $\lambda_* N^i$, and then we apply
$\sigma$. 
\end{enumerate}
\end{proof}

\subsection{Products} \label{subsec:pro}

In this section, we want to define products on $AD^{n}_{\tau_0}
(X)$ and $D^{n}_{\tau_0} (X)$, compatibly with the products in
the Chow groups and in the Deligne cohomology. To this aim we
follow the pattern explained in \cite{B} and \cite{EV}. 

\begin{defn} \label{dfn:pro}

Let $\alpha \in \R$. We define 
$$
C(m)_{\tau_0} \times C (n)_{\tau_0} \> \cup_{\alpha} >> C (m
+n)_{\tau_0} 
$$
by 
$$
\begin{array}{lll}
x \cup_{\alpha} y & = \{ x,y\} & x \in \sK_m , y \in \sK_n \\
&& \\
& = 0 & x \in \sK_m, y \in N^{\geq 2n} [n] \\
&& \\
& = (1 - \alpha) d \log x \wedge y & x \in \sK_m, y \in N^{\geq
n} [n] \\
&& \\ 
& =0 & x \in N^{\geq 2m} [m], y \in \sK_n \\
&& \\
& = x \wedge y & x \in N^{\geq 2m} [m] , y \in N^{\geq 2n} [n]
\\ 
&& \\
& = (-1)^{\deg x} \alpha x \wedge y & x \in N^{\geq 2m} [m] , y
\in N^{\geq n} [n] \\
&& \\
& = \alpha x \wedge d \log y & x \in N^{\geq m} [m], y \in \sK_n
\\
&& \\
& = (1- \alpha) x \wedge y & x \in N^{\geq m} [m], y \in N^{\geq
2n} [n] \\
&& \\
& = 0 & x \in N^{\geq m} [m], y \in N^{\geq n} [n] 
\end{array} 
$$
\end{defn} 

\begin{prop}
\begin{enumerate}
\item[1)] These formulae define for each $\alpha \in \R$ a
product, compatibly with the product on $\sK_n, N^{\geq 2n}$,
and with the product 
\begin{gather}
\sK_m \times N^{\geq 2n} \> d \log \times 1 >> \Omega^{\geq
m}_{X} [m] \times N^{\geq 2n} \> \tau_0 >> N^{\geq m + 2n} [m]
\notag \\
N^{\geq 2m} \times \sK_n \> 1 \times d \log >> N^{\geq 2m}
\times \Omega^{\geq n} [n] \> \tau_0 >> N^{2m +n} [n] \notag
\end{gather}
\item[2)] One has $x \cup_{\alpha} y = (-1)^{\deg x \deg y} y
\cup_{(1- \alpha)} x$. 
\item[3)] For $\alpha $ and $\beta$, the products
$\cup_{\alpha}$ and $\cup_{\beta}$ are homotopic. We denote by
$\cup$ the induced product in cohomology. In particular
\begin{enumerate}
\item[a)] $\cup$ is commutative on $AD^{n}_{\tau_0} (X)$. 
\item[b)] The restriction of $\cup$ to $\H^n (X, N^{\infty}
\sK_n)$ (\ref{prop:diag}, 2 ii)) is given by
$$
\begin{array}{ll}
\{ x,y\} & x \in \sK_m, y \in \sK_n \\
& \\
d \log x \wedge y & x \in \sK_m , y \in N^{\geq n} [n] \\
& \\
0 & \mbox{otherwise} 
\end{array} 
$$
\end{enumerate}
\end{enumerate}
\end{prop}

\begin{proof} The verification is exactly as in \cite{EV},
(3.2), (3.3), (3.5), where one replaces 
$$
\begin{array}{lll}
\Z (n) & \mbox{by} & \sK_n \\
& & \\
F^n & \mbox{by} & N^{\geq 2n} [n] \\
& & \\
\Omega^{\bullet}_{X} & \mbox{by} & N^{\geq n} [n]. 
\end{array} 
$$
In particular, the homotopy between $\cup_{\alpha}$ and
$\cup_{\beta}$ is given by
$$ 
\begin{array}{lll}
h (x \otimes y) &=  (-1)^{\mu} (\alpha - \beta) x \wedge y 
& \mbox{if} \ \ x \in (N^{\geq m} [m])^{\mu-1} y \in (N^{\geq n}
[n])^{\mu' -1} \\
&= 0 &\ \ \ \mbox{otherwise} 
\end{array}
$$
for
\begin{equation*}
h: (C (m)_{\tau_0} \otimes_{\Z} C (n)_{\tau_0})^{\ell} \to C (m
+ n)^{\ell}_{\tau_0}. 
\end{equation*}
\end{proof}

\section{The splitting principle for bundles with connection}
\label{sec:spl}  

\subsection{The $\tau$ complex} \label{subsec:tau}

Let $\tau_0 : \Omega^{\bullet}_{X} (\log D) \to N^{\bullet}$ and
$(E, \nabla)$ be as in \ref{sec:not} 4), 6), with $B^b = B^a
\wedge \Omega^{b - a}_{X} (\log D)$. Let $\pi: P : = \P (E)
\to X$ be the projective bundle of $E$, and $D' : = \pi^{-1}
(D)$. 

\begin{defn} \label{dfn:tauform}
$$
\Omega^{n}_{\P, \tau_0} (\log D') = \frac{\pi^* N^n \oplus
\Omega^{n}_{\P} (\log D') / \pi^* B^a \wedge \Omega^{n-a}_{\P}
(\log D')}{\pi^* (\Omega^{n}_{X} (\log D)/ B^a \wedge
\Omega^{n-a}_{X} (\log D))} . 
$$
\end{defn}

For example, if $\Omega^{\bullet}_{X} (\log D)$ surjects onto
$N^{\bullet}$, then this is 
$$
\sF^n = \Omega^{n}_{\P} (\log D')/ \pi^* B^a \wedge
\Omega^{n-a}_{\P} (\log D'),
$$
and in general it is the push-down of this sheaf via
$$
\begin{CD}
\pi^* (\Omega^{n}_{X} (\log D) / B^a \wedge \Omega^{n-a}_{X}
(\log D)) \> \subset >> \sF^n \\
\V \subset VV \noarr \\
N^n. \noarr 
\end{CD}
$$
The Leibniz rule implies that $\pi^* B^a \wedge
\Omega^{n-a}_{\P} (\log D')$ is a subcomplex of
$\Omega^{\bullet}_{\P} (\log D')$, and therefore $\tau_0$
induces a map of differential graded algebras, still denoted by
$\tau_0$: 

\begin{defn} 
$$
\tau_0 : (\Omega^{\bullet}_{\P} (\log D') , d) \to
(\Omega^{\bullet}_{\P, \tau_0} (\log D'), \tau_0 \circ d).
$$
\end{defn}

Then $\nabla$ induces a $\tau_0$ connection, still denoted by
$\nabla$: $$
\nabla: \pi^* E \to \Omega^{1}_{\P , \tau_0} (\log D') \otimes
\pi^* E.
$$
Recall from \cite{EI}, \S \ 2, that the connection 
$$
\nabla : E \to N^1 \otimes E
$$
induces a splitting 
$$
\tau : \Omega^{1}_{\P, \tau_0} (\log D') \to \pi^* N^1
$$
of the exact sequence
$$
0 \to \pi^* N^1 \to \Omega^{1}_{\P, \tau_0} (\log D') \to
\Omega^{1}_{\P /X} \to 0, 
$$
such that the $\pi^* N^1$ valued connection $\tau \circ \nabla$
on $\pi^* E$ respects the canonical filtration 
$$
0 \to \Omega^{1}_{\P /X} (1) \to \pi^* E \to \sO (1) \to 0, 
$$
i.e.: 
$$
\tau \circ \nabla (\Omega^{1}_{\P / X} (1)) \subset \pi^* N^1
\otimes \Omega^{1}_{\P /X} (1). 
$$
The induced connection on $\Omega^{1}_{\P /X} =
\Omega^{1}_{\P/X} (1) \otimes \sO (1)^*$ is then given by
applying first the splitting into $\Omega^{1}_{\P, \tau_0} (\log
D')$, then applying the differential $d$ of
$\Omega^{\bullet}_{\P , \tau_0} (\log D')$, then projecting onto
the factor $\pi^* N^1 \otimes \Omega^{1}_{\P/X}$ of
$\Omega^{1}_{\P, \tau_0} (\log D')$. We write for short

\begin{defn} 
$$
\tau d: \Omega^{1}_{\P /X} \to \pi^* N^1 \otimes \Omega^{1}_{\P
/ X}.
$$
\end{defn}

\begin{lem} \label{lemma:wt}

In the splitting 
$$
\Omega^{n}_{\P , \tau_0} (\log D') = \oplus^{n}_{a=0}
\Omega^{a}_{\P/X} \otimes \pi^* N^{n-a}
$$
one has 
\begin{multline*}
\tau \circ d (\Omega^{a}_{\P/X} \otimes \pi^* N^{n-a}) \\
\subset \Omega^{a-1}_{\P /X} \otimes \pi^* N^{n-a+2} \oplus
\Omega^{a}_{\P /X} \otimes \pi^* N^{n-a+1} \\
\oplus \Omega^{a+1}_{\P / X} \otimes \pi^* N^{n-a}
\end{multline*}
\end{lem} 

\begin{proof}
The map $\tau \circ d$ has 3 components: 
\begin{gather}
d_{\text{rel}} : \Omega^{1}_{\P/X} \to \Omega^{2}_{\P/X} \notag \\
\tau d \notag \\
\mu : \Omega^{1}_{\P/X} \to \pi^* N^2, \notag
\end{gather}
where $\mu$ is $\sO_{\P}$ linear (see \cite{BE}, (4.3.2) for the
study of $\mu$). One just applies the Leibniz rule. 
\end{proof}

\begin{defn} \label{dfn:xi}

\begin{enumerate}
\item[1)] We define a decreasing filtration $\Phi$ on
$\Omega^{n}_{\P, \tau_0} (\log D')$ by 
$$
\Phi^{\ell} \Omega^{n}_{\P , \tau_0} (\log D') =
\oplus^{n}_{a=\ell} \Omega^{a}_{\P/X} \otimes \pi^* N^{n-a}
$$
for $\ell \leq n$
\begin{gather}
\Phi^{\ell} \Omega^{n}_{\P, \tau_0} (\log D') = 0 \ \
\mbox{for} \ \ \ell > n, \notag \\
\Phi^0 \Omega^{n}_{\P, \tau_0} (\log D') = \Omega^{n}_{\P,
\tau_0} (\log D') \notag
\end{gather}
fulfilling 
$$
d \Phi^{\ell} \Omega^{n}_{\P, \tau_0} (\log D') \subset
\Phi^{\ell -1} \Omega^{n+1}_{\P, \tau_0} (\log D')
$$
\item[2)] We define the complex 
$$
\begin{array}{lll}
(M^{\geq n})^{\ell} & = 0 & \ell < n \\
& \\
& = \Omega^{\ell}_{\P , \tau_0} (\log D') / \Phi^{2n- \ell}
\Omega^{\ell}_{\P , \tau_0} (\log D')  & \\
& \\
&( = \oplus_{a < 2n - \ell} \Omega^{a}_{\P /X} \otimes \pi^*
N^{\ell -a}) &.
\end{array}
$$
In particular, $(M^{\geq n})^{\ell} =0$ for $\ell \geq 2n$. 
\item[3)] We define the complex 
$$
\begin{array}{ll}
\Gamma (n)_{\tau} & = {\rm cone} \ (\sK_n \> \tau \circ d \log
>> M^{\geq n} [n] ) [-1] \\
& \\
& = {\rm cone} \ (\sK_n \oplus (M^{\geq n})^{\geq 2n} [n] \>
\tau \circ d \log \oplus 0 >> M^{\geq n} [n]) [-1]
\end{array}
$$
endowed with the map 
$$
\pi^{-1} : C (n)_{\tau_0} \to R \pi_* \Gamma (n)_{\tau}
$$
induced by 
\begin{gather}
\pi^{-1} \sK_n \to \sK_n \notag \\ 
\pi^{-1} N^{\geq 2n} \to 0 = (M^{\geq n})^{\geq 2n}, \notag \\
\pi^{-1} N^{\geq n} \to M^{\geq n}. \notag
\end{gather}
\item[4)] We define the multiplication 
$$
\Gamma (m)_{\tau} \times \Gamma (m)_{\tau} \to \Gamma (m
+n)_{\tau} 
$$
as in \ref{dfn:pro}, replacing $N^{\geq n}$ by $M^{\geq n}$,
and $N^{\geq 2n}$ by $0$, and observing that for $x \in 
\Omega^{m}_{\P, \tau_0} (\log D')$, $y \in \oplus_{a < 2 n-\ell}
\Omega^{a}_{\P /X} \otimes \pi^* N^{\ell - a}$, then 
\begin{multline*}
x \wedge y \in \\
\oplus_{a < 2 n + m - \ell} \Omega^{a}_{\P/X}
\otimes \pi^* N^{m + \ell -a} \\
= \oplus_{a < 2 (n + m) - (\ell + m)} \Omega^{a}_{\P /X} \otimes
\pi^* N^{m + \ell -a} .
\end{multline*}
Thus concretely 
$$
\begin{array}{lll}
x \cup_{\alpha} y &  = \{ x,y\} \ \ \mbox{for} \ \ x \in \sK_m, y
\in \sK_n \\
& \\
& = (1- \alpha) d \log x \wedge y \in (M^{\geq m +n})^{m + \ell}
\\ 
& \\
& \hspace{1cm} \mbox{for} \ \ x \in \sK_m, y \in (M^{\geq
n})^{\ell} \\
& \\
& = \alpha x \wedge d \log y \in (M^{\geq m +n})^{n + \ell} \\
& \\
& \hspace{1cm} \mbox{for} \ \ x \in (M^{\geq m})^{\ell} , y \in
\sK_n \\
& \\
& = 0 \ \mbox{for} \ \ x \in M^{\geq m} [m] , y \in M^{\geq n}
[n]. 
\end{array} 
$$
Again $\cup_{\alpha}$ does not depend on $\alpha$ on cohomology.
We denote it by $\cup$. In particular, the product induces an
action
\begin{multline*} 
\H^a (\P, M^{\geq b} [b]) \times \H^{a'} (\P, \Gamma
(b')_{\tau}) \\
\> \cup >> \H^{a+a'} (\P, M^{\geq (b+b')} [b+b'])
\end{multline*}
making Image $\H^{\bullet -1} (\P, N^{\geq \bullet \bullet}
[\bullet \bullet])$ in $\H^{\bullet} (\P, \Gamma (\bullet
\bullet)_{\tau})$ into an ideal of square $0$. 
\item[5)] The group $\H^1 (\P, \Gamma (1)_{\tau}) = \H^1 (\P,
\sK_1 \to \pi^* N^1)$ is the group of isomorphism classes of
rank 1 bundles with a $\pi^* N^1$ valued connection. We denote
by $\xi$ the class of $(\sO (1), \tau \circ \nabla)$ in it
(\ref{subsec:tau}) 
\end{enumerate}
\end{defn}

\begin{thm} \label{thm:splitting}

The maps 
$$
\pi^{-1} : AD^{a}_{\tau_0} (X) \to \H^a (\P, \Gamma (a)_{\tau})
$$
are injective for all $a \geq 0$. One has a splitting 
$$
\H^m (\P, \Gamma (n)_{\tau}) = \oplus^{r-1}_{i=0} \H^{m-i} (X, C
(n-i)_{\tau_0} ) \cup \xi^i
$$
(where $r$ is the rank of $E$ (\ref{sec:not}, 6)). 
\end{thm}

\begin{proof} 

We denote by $[\xi]$ the class of $\sO (1)$ in $H^1 (\P ,
\sK_1)$. One has the projective bundle formula for 
the cohomology of the sheaves
$\sK$ 
$$
\H^m (\P, \sK_n) = \oplus^{r-1}_{i=0} H^{m-i} (X, \sK_{n-i} )
\cup [\xi ]^i.
$$
One also has obviously
$$
\H^m (\P, M^{\geq n} [n]) = \oplus^{r-1}_{i=0} \H^{m-i} (X,
N^{\geq n-i} [n-i] ) \cup \xi^i .
$$
One regards the exact sequences 
\begin{multline*}
\to H^{m-1} (\P , \sK_n) \to \H^{m-1} (\P, M^{\geq n} [n]) \to
\\
\H^m (\P, \Gamma (n)_{\tau}) \to H^m (\P, \sK_n ) \to \H^m (\P,
M^{\geq n} [n]) \to  
\end{multline*}
and
\begin{multline*}
\to \oplus^{r-1}_{i=0} H^{m-1-i} (X, \sK_{n-i}) \to
\oplus^{r-1}_{i=0} \H^{m-1-i} (X, N^{\geq r -1 -i} [r-1-i])
\\
\to \oplus^{r-1}_{i=0} \H^{m-i} (X, C (n-i)_{\tau_0}) \to
\oplus^{r-1}_{i=0} H^{m-i} (X, \sK_{n-i})  \\
\to \oplus^{r-1}_{i=0} \H^{m-i} (X, M^{\geq (n-i)} [n-i]) \to . 
\end{multline*}
The map $\pi^{-1}$ induces an isomorphism on the $\sK$ and $M$
cohomology, thus on the $\Gamma$ cohomology as well. 
\end{proof} 

\begin{defn} \label{dfn:class}

Let $(E, \nabla), \tau_0, \xi$ be as in \ref{dfn:xi}. We define
$$c_0 (E, \nabla) = 1 \in AD^{0}_{\tau_0} (X) = \Z$$ 
and $$c_n (E,
\nabla) \in AD^n_{\tau_0}(X), n = 1, \ldots , r$$ by the formula 
$$
\sum^{r}_{n=0} (-1)^n c_n (E, \nabla) \cup \xi^{r-n} =0
$$
via \ref{thm:splitting}. We define $c_n (E, \nabla) =0$ for $n > m$. 
\end{defn}

\begin{thm} \label{thm:char} \ \\
\begin{enumerate}
\item[1)] {\bf Functoriality:} Let $\pi : Y \to X$ be a morphism of
smooth varieties, such that $D' = \pi^{-1} D$ is a normal
crossing divisor. Let $\tau_0: \Omega^{\bullet}_{Y} (\log D')
\to \Omega^{\bullet}_{Y, \tau_0} (\log D')$ be the map of
differential graded algebras with $\Omega^{n}_{Y, \tau_0} (\log
D')$ defined as in \ref{dfn:tauform} by 
$$
\Omega^{n}_{Y, \tau_0} (\log D') = \frac{\pi^* N^n \oplus
\Omega^{n}_{Y} (\log D')/ \pi^* B^a \wedge \Omega^{n-a}_{Y}
(\log D')}{\pi^* (\Omega^{n}_{X} (\log D) / B^a \wedge
\Omega^{n-a}_{X} (\log D))}. 
$$
Then $\nabla$ induces on $\pi^* E$ a $\Omega^{1}_{Y, \tau_0}
(\log D')$ valued connection $\pi^* \nabla$, and $\pi^{-1}$
induces a map
$$
\pi^{-1} : AD^{n}_{\tau_0} (X) \to AD^{r}_{\tau_0} (Y).
$$
Then $c_n (\pi^* E, \pi^* \nabla) = \pi^{-1} c_n (E, \nabla)$.
For all further maps 
$$
\tau : \Omega^{\bullet}_{Y, \tau_0} (\log D') \to M^{\bullet}
$$
of differential graded algebras, one has 
$$
c_n (\pi^* E, \pi^* \nabla) = {\rm Im} \ c_n (\pi^* E, \pi^*
\nabla) 
$$
in $AD^{n}_{\tau \circ \tau_0} (X)$. 
\item[2)] {\bf First class:} If $E$ is of rank 1, and $\tau_0 ,
\nabla$ are as in \ref{sec:not}, 4), 6), then $c_1 (E, \nabla)$
is the isomorphism class of $(E, \nabla)$ in the group of
isomorphism classes of rank 1 bundles with a $\tau_0$ connection
$AD^{1}_{\tau_0} (X) = \H^1 (X, \sK_1 \to \pi^* N^1 )$. 
\item[3)] {\bf Whitney product formula:} Let 
$$
0 \to (E'', \nabla'') \to (E, \nabla) \to (E', \nabla') \to 0
$$
be an exact sequence of bundles with $\tau_0$ connection, that
is the bundle sequence is exact and 
$$
\nabla E'' \subset N^1 \otimes E'' \ , \ \ \nabla'' = \nabla|_{E''}
, \nabla'
$$
is the quotient connection. Then one has 
$$
c_n (E, \nabla) = \oplus^{n}_{a=0} c_a (E', \nabla') \cup
c_{n-a} (E'', \nabla'') 
$$
in $AD^{n}_{\tau_0} (X)$ for all $n \geq 0$. 
\item[4)] {\bf Characterization of the classes:} 
The classes $c_n (E, \nabla) \in AD^{n}_{\tau_0} (X)$
are uniquely determined by the functoriality property, the
definition of the first class and the Whitney product formula. 
\end{enumerate}
\end{thm}

\begin{proof} 
1) and 2) are clear, and 4) is clear once one knows 3). For 3) ,
we mimic the classical proof \cite{I}, 6.10. One has 
$$
\begin{CD}
\P' : = \P (E') \> i >> P \< j << U: \P - \P' \\
\noarr \noarr \V V p V \\
\noarr \noarr \P'' : = \P (E'') 
\end{CD}
$$
where $p$ is an affine filtration with the obvious notations
$\Gamma (m)_{\tau'}, r'' , ...$ one defines 
\begin{gather}
\alpha = \sum^{r'}_{n=0} (-1)^n c_n (E', \nabla') \cup
\xi^{r'-n} \in \H^{r'} (\P, \Gamma (r')_{\tau} ) \notag \\
\beta = \sum^{r''}_{n=0} (-1)^n c_n (E'', \nabla'') \cup
\xi^{r''-n} \in \H^{r''} (\P, \Gamma (r'')_{\tau}) \notag
\end{gather}
Then $\gamma = \alpha \cup \beta \in \H^r (\P, \Gamma
(r)_{\tau})$ fulfills 
$$
j^* \beta = p^{-1} (\sum^{r''}_{n=0} (-1)^n c_n (E'', \nabla'')
\cup \xi^{\prime \prime^{r'' -n}} ) =0
$$
and $\alpha|_{\P'} =0$. Let $\beta' \in \H^{r}_{\P'} (\P, \Gamma
(r)_{\tau})$, with $\rho (\beta') = \beta$ in the localization
sequence 
$$
\H^{r}_{\P'} (\P, \Gamma (r)_{\tau}) \> \rho >> \H^r (\P, \Gamma
(r)_{\tau}) \> j^* >> \H^r (U, \Gamma(r)_{\tau}).
$$
Since the product 
$$
\H^a (\P, \Gamma (b)_{\tau}) \times \H^{a'} (\P, \Gamma
(b')_{\tau}) \to \H^{a+a'} (\P, \Gamma (b + b')_{\tau})
$$
induces a compatible product
$$
\H^a (\P, \Gamma (b)_{\tau}) \times \H^{a'}_{\P'} (\P, \Gamma
(b')_{\tau}) \to \H^{a + a'}_{\P'} (\P, \Gamma (b+b')_{\tau})
$$
one has $\gamma = \rho (\alpha \cup \beta')$. On the other hand,
one has a Gysin isomorphism
$$
\H^{a-r''} (\P', \Gamma (b- r'')_{\tau''}) \> i_* >>
\H^{a}_{\P'} (\P, \Gamma (b)_{\tau}).
$$
To see this, one observes that $i_*$ induces an isomorphism on
the $\sK$ cohomology. For the $M^{\geq b}$ cohomology, one makes
a d\'evissage with the stupid filtration. This then reduces to
the obvious isomorphism 
$$
\H^a (\P', \Omega^{b}_{\P'/X} \otimes \pi^{\prime *} N^c) \> i_* >>
\H^{a+r''}_{\P'} (\P, \Omega^{b+r''}_{\P/X} \otimes \pi^* N^c ).
$$
We write $\beta' = i_* (\beta'')$ with $\beta'' \in \H^{r'} (\P,
\Gamma (r')_{\tau})$, and $\gamma = \rho (\alpha \cup i_*
\beta'')$. It remains to see that 
$$
x \cup i_* y = i_* (x |_{\P'} \cup y)
$$
for $x \in \H^a (\P, \Gamma (b)_{\tau})$, $y \in \H^{a'-r''}
(\P', \Gamma (b' - r'')_{\tau'})$, as $\alpha |_{\P'} =0$. One
checks that the diagram 
$$
\begin{CD}
\Gamma (b)_{\tau} \times {\rm cone} \ (\Gamma (b')_{\tau} \to R
j_* \Gamma (b')_{\tau} ) [-1] \> \cup >> {\rm cone} \ (\Gamma
(b+b') \to R j_* \Gamma (b+b')_{\tau}) [-1] \\
\V \mbox{restriction}\times \mbox{residue} VV \V V \mbox{residue} V \\
\Gamma (b)_{\tau'}|{\P'} \times \Gamma (b' - r'')_{\tau'}|{\P'} \> \cup >>
\Gamma (b+b' -r'')_{\tau'}|{\P'}
\end{CD}
$$
is commutative. 
\end{proof}

\subsection{Comparison with the algebraic Chern-Simons
invariants of \cite{BE}} 

\begin{thm} \label{thm:CS}
\begin{enumerate}
\item[1)] Let $(\tau_0 , E, \nabla)$ be as in
\ref{thm:splitting}, with $\tau_0$ surjective. Then the image of
$c_n (E, \nabla)$ in $H^0 (X, N^{2n}_{cl})$ (\ref{prop:diag}
2,i)) is $P_n (\nabla^2 , \ldots , \nabla^2)$ (\ref{sec:not},
10)). 
\item[2)] Assume $\tau_0 =$ identity. Then the image of $c_n (E,
\nabla)$ under
$$
AD^{n}_{\tau_0} (X) \to H^0 (X, \Omega^{2n-1}_{X} (\log D)/d
\Omega^{2n-2}_{X} (\log D))
$$
coincides with the algebraic Chern-Simons class $w_n (E,
\nabla)$ defined in \cite{BE}, \S \ 2. 
\end{enumerate}
\end{thm}

\begin{proof} 
\begin{enumerate}
\item[1)] Let $p: G \to X$ be the flag bundle of $E, \Delta =
p^{-1} (D)$. Define as in \cite{BE}, (4.3.3) the sheaf 
$$
\bar{M}^{2n-\ell} = p^* N^{2n+\ell} / \mu (\Omega^{1}_{G/X}
\otimes p^* N^{2n+\ell -2})
$$
where $\mu$ is as in the proof of \ref{lemma:wt}, for $\ell \geq
0$. One defines $\Omega^{\bullet}_{G, \tau_0} (\log D)$ as in
\ref{thm:char}, 1), and $\Gamma (n)_{\tau}$ similarly. Let 
$$
\bar{\Gamma} (n)_{\tau} = {\rm cone} \ (\sK_n \to \bar{M}^{\geq
n} [n]) [-1] $$
with
$$
\begin{array}{lll}  
 (\bar{M}^{\geq n})^{\ell} & = (M^{\geq n})^{\ell} & \ell \leq 2n-1 \\
  & = \bar{M}^{\ell} &  \ell \geq 2n. 
\end{array}
$$
Then $\bar{\Gamma} (n)_{\tau}$ maps to $\Gamma (n)_{\tau}$. We
lift the product in $\Gamma (n)_{\tau}$ to $\bar{\Gamma}
(n)_{\tau}$ defined in \ref{dfn:xi} by setting 
$$
\begin{array}{lll}
x \cup_{\alpha} y & = \{ x,y\} & x \in \sK_m , y \in \sK_n \\
& & \\
& = (1-\alpha) d \log x \wedge y & x \in \sK_m, y \in
(\bar{M}^{\geq n} [n])^{\ell} \\
& & \\
& = \alpha x \wedge d \log y & x \in (\bar{M}^{\geq m}
[m])^{\ell},  y
\in \sK_n \\
& & \\
& = 0 & x \in \bar{M}^{\geq m} [m], y \in \bar{M}^{\geq n} [n] .
\end{array} 
$$
Here the product $d \log x \wedge y$ $(x \wedge d \log y)$ means
$\tau (d \log x \wedge y)$ $(\tau (x \wedge d \log y))$ for $m+
\ell \geq 2 (m+n)$ $(\ell + n \geq 2 (m+n))$. The product being
compatible with the natural product of $(\bar{M}^{\geq n})^{\geq
2n}$, the image of $p^{-1} c_n (E, \nabla)$ in $H^0 (G,
\bar{M}^{2n} \to \bar{M}^{2n+1} \to ...)$ is the $n$-th
symmetric product of the image $F(\ell_i)$ of 
$$
\ell_i = c_1 (L_i, \tau \circ \nabla) \in \H^1 (G, \sK_1 \to p^* N^1
) 
$$
in 
$$
\H^2 (G, (\bar{M}^{\geq 1})^{\geq 2} ) = \H^0 (G \frac{p^*
N^2}{\mu (\Omega^{1}_{G/X})} \to ...). 
$$
In other words, considering as in \cite{BE}, \S \ 4, the
quotient differential graded algebra 
$$
\Omega^{\bullet}_{G, \tau_0} (\log \Delta) \to (\sO_{G} \to
\bar{M}^{\geq 1}),
$$
$F(\ell_i)$ is the curvature of the $(\bar{M}^{\geq 1})^1 = p^*
N^1$ valued connection $\tau \circ \nabla$ on $L_i$. As in
\cite{BE}, (4.7.8), this is exactly $p^* P_n (\nabla^2 , \ldots
, \nabla^2)$. Now, as $N^{2n}$ si locally free, $H^0 (X,
N^{2n})$ injects into $H^0 (U, N^{2n})$ for any open $\phi \neq
U \subset X$ . Thus we may assume that $\nabla : E \to N^1
\otimes E$ lifts to $\bar{\nabla} : E \to \Omega^{1}_{X} \otimes
E$ as $\Omega^{1}_{X} (\log D) \to N^1$ is surjective, and it is
enough to prove 1) for $\bar{\nabla}$. As in \cite{BE}, \S \ 4,
we may further assume that there is a morphism $\varphi : X \to
T$, where $T$ is an affine space, such that $(E, \bar{\nabla}) =
\varphi^* (\sE , \psi)$, where $\sE$ is the trivial bundle, and
$\psi$ is such that if $q: P \to T$ is the flag bundle of $\sE$,
and $(\bar{M}^{\geq 1})$ is defined on $p$, then the map
$\Omega^{2n}_{T} \to p_* (\bar{M}^{\geq 1})^{2n}$ is injective
(\cite{BE}, Proposition 4.9.1). Since 
\begin{multline*}
AD^n (T) = \frac{H^0 (T, \Omega^{2n-1}_{T})}{d H^0 (T,
\Omega^{2n-2}_{T})} = H^0 (X, \frac{\Omega^{2n-1}_{T}}{d
\Omega^{2n-2}_{T}}) \\
\subset H^0 (T, \Omega^{2n}_{T}) \subset H^0 (P,
(\bar{M}^{\geq 1})^{2n})
\end{multline*}
we see that 
$$
{\rm Im} \ c_n (\sE, \psi) \ \mbox{in} \ H^0 (P,
(\bar{M}^{\geq 1})^{2n}) 
= p^* P_n (F (\psi), \ldots , F (\psi) )
$$
implies that 
$$
{\rm Im} \ c_n (\sE, \psi) \ \mbox{in} \ H^0 (T,
\Omega^{2n}_{T}) \\
= P_n (F (\psi), \ldots , F(\psi)).
$$
\item[2)] From the exact sequence 
\begin{multline*}
0 \to H^0 (X -D, \sH^{2n-1}) \to H^0 (X, \frac{\Omega^{2n-1}_{X}
(\log D)}{d \Omega^{2n-1}_{X} (\log D)} ) \\
\to H^0 (X, \Omega^{2n}_{X} (\log D))
\end{multline*}
and the injectivity 
$$
H^0 (X - D, \sH^{2n-1}) \to H^0 (U - D, \sH^{2n-1})
$$
for any open $0 \neq U \subset X$ (\cite{BO}, (4.2.2)) one sees
that the question is local. We may check it on $(T, \sE, \psi)$
as in 1). But on $T$, the algebraic Chern-Simons invariant is
the only class mapping to $P_n (F(\psi), \ldots , F(\psi))$. 
\end{enumerate}
\end{proof}

\subsection{Comparison with classes of bundles with flat
connections} 

\begin{thm} \label{thm:flat}
Let $(E, \nabla, \tau_0)$ be as in \ref{thm:splitting}, with
$\tau_0$ surjective and with $\nabla^2 =0$. Then the class $c_n
(E, \nabla) \in \H^n (X, N^{\infty} \sK_n)$ (\ref{prop:diag})
are the same as the classes defined in \cite{EII} (see also
\cite{BE}, (3.11)). 
\end{thm}

\begin{proof} 
We consider the flag bundle $p: G \to X$, with $\Delta = p^{-1}
(D)$. We have the quotient maps 
$$
\Omega^{\bullet}_{G} (\log \Delta) \to \Omega^{\bullet}_{G,
\tau_0} (\log \Delta) \to p^* N^{\bullet} 
$$
when $\Omega^{\bullet}_{G, \tau_0} (\log \Delta)$ is as in
\ref{thm:char}, 1). We also have the filtration $\Phi$ on
$\Omega^{\bullet}_{G, \tau_0} (\log \Delta)$ defined as in
\ref{dfn:xi}, for $G$ replacing $\P$, and the corresponding
complex $\Gamma (n)_{\tau}$. One has a quotient map 
$$
\Gamma (n)_{\tau} \to {\rm cone} \ (\sK_n \oplus (p^* N^{\geq
2n}) [n] \to (p^* N^{\geq n}) [n]) [-1]
$$
and the multiplication $\cup_0$ on $\Gamma (n)_{\tau}$ defined in
\ref{dfn:xi} maps to 
$$
\begin{array}{lll}
x \cup_0 y & = \{ x, y\} & x \in \sK_m, y \in \sK_n \\
& & \\
& = \tau d \log x \wedge y & x \in \sK_m , y \in p^* N^{\ell} \\
& & \\
& = 0 & \mbox{otherwise} .
\end{array} 
$$
By the Whitney product formula \ref{thm:char}, 3), the class
$c_n (E, \nabla) \in AD^{n}_{\tau_0} (X)$ is the $n$-th
symmetric function of the rank one subquotients 
$$(L_i, \tau
\circ \nabla) 
\in \H^1 (G, \sK_1 \to (p^* N^1)_{cl}) 
\subset \H^1(G, \sK_1 \to p^* N^1).
$$
Since the multiplication on $\Gamma (n)_{\tau}$, restricted to
the elements of $\H^n (G,\Gamma (n)_{\tau})$ mapping to $\H^n
(G, \sK_n \to (p^* N^n)_{cl})$, maps to the multiplication on
$(\sK_n \to (p^* N^n)_{cl})$ defined in \cite{EII}, p. 51 in
order to construct the classes of flat bundles, one has that the
image of $c_n (E, \nabla)$ in $\H^n (G, p^* N^{\infty} \sK_n)$
is the class defined in \cite{EII}. One concludes by the
commutative diagram 
$$
\begin{CD}
\H^n (X, N^{\infty} \sK_n) \> \subset >> AD^{n}_{\tau_0} (X) \\
\V \subset VV \V \subset VV \\
\H^n (G, \sK_n \to p^*N^n \to \cdots \to p^*N^{2n-1}) \< <<  
\H^n (G, \Gamma
(n)_{\tau}), 
\end{CD}
$$
where the left vertical injective arrow is 
the composition of the injection
$\H^n(X, N^{\infty} \sK_n) \to \H^n(G, p^*N^{\infty} \sK_n)$ 
as in \cite{EII},
p.51, p.52, with the injection
$\H^n(G, p^*N^{\infty} \sK_n) \to
\H^n (G, \sK_n \to p^*N^n \to \cdots \to p^*N^{2n-1}), $
knowing that
$$ \rm{Im} \ c_n(E, \nabla) \ \rm{in} \ H^0(X, N^{2n}_{cl})=
P_n (\nabla^2 , \ldots , \nabla^2) =0
$$
\ref{thm:CS}, 1). 
\end{proof}\section{The Weil Algebra}  

\subsection{Construction of Beilinson-Kazhdan} \label{subsection:BK}

\begin{enumerate}
\item[1.] In the unpublished note \cite{BK}, Beilinson and
Kazhdan construct a complex generalizing the Weil homomorphism
\begin{gather}
(S^n \sG^* )^G \to H^{2n}_{DR} (X) \notag \\
P \mapsto P (\nabla^2 , \ldots , \nabla^2) \notag
\end{gather}
associated to $(E, \nabla)$. We give a short account of their
construction. Instead of considering reductive groups in
general, we consider the subgroup $G_e \subset G$  of matrices fixing the
flag 
$$
e: 0 \to k^{r''} \to k^r \to k^{r'} \to 0.
$$
We denote by $(E, e)$ and exact sequence 
$$
0 \to E'' \to E \to E' \to 0
$$ 
with rank $E'' = r''$, with the convention $(E, 0) = E$,
$G_0=G$ (see
\ref{sec:not}, 10), 11)). Let 
$$
0 \to \Omega^{1}_{X} \> \iota >> \Omega^{1}_{X, (E, e)} \> \pi
>> \sG^{*}_{(E, e)} \to 0 
$$
be the (functorial) Atiyah sequence of $(E, e)$, whose splitting
is equivalent to a $\Omega^{1}_{X}$ valued connection on $E$
compatible with $e$. Then $\Omega^{1}_{X, (E,e)} = (p_*
\Omega^{1}_{\sE_e})^G_e$, where $p: \sE_e \to X$ is the total
space of the $G_e$ torseur (\ref{sec:not}, 10)) associated to $(E,
e)$. Let $\Omega^{\bullet}_{X, (E, e)}$ be the sheaf of
commutative differential graded algebras generated by $\sO_X$ in
degree 0, $\Omega^{1}_{X}$ in degree 1, such that for $f \in
\sO_X, df$ is the K\"ahler differential in $\Omega^{1}_{X}
\subset \Omega^{1}_{X, (E, e)}$. Then 
$$
\Omega^{n}_{X,(E,e)} = \oplus_{a +b =n} \Omega^{a,b}_{X, (E, e)} 
$$
with
$$
\Omega^{a,b}_{X, (E,e)} = \Lambda^{a-b} \Omega^{1}_{X, (E,e)}
\otimes S^b \sG^{*}_{(E,e)} 
$$
with differential $d = d' + d''$,
$$
 d' : \Omega^{a,b}_{X, (E,e)} \to \Omega^{a+1,
b}_{X, (E,e)}
$$
being induced from $\sE_e$ via $\Lambda^i \Omega^{1}_{X, (E,e)} = 
(p_* \Omega^{i}_{\sE_e} )^G$, and $$d'': \Omega^{a,b}_{X, (E,e)} 
\to \Omega^{a,b+1}_{X,(E,e)}$$
being defined by  $ d'' | \Omega^{10}_{X, (E,e)} =
d'' | \Omega^{1}_{X,(E,e)} = \pi$, 
$$
d'' (x_1 \wedge \ldots \wedge x_{a-b} \otimes \sigma) =
\sum^{a-b}_{i=1} (-1)^i x_1 \wedge \ldots \wedge \hat{x_i}
\ldots \wedge x_{a-b} \pi (x_i) \cdot \sigma
$$
for
$$
x_1 \wedge \ldots \wedge x_{a-b} \in \Lambda^{a-b}
\Omega^{1}_{X, (E,e)}, \sigma \in S^b \sG^{*}_{(E,e)} .
$$
\item[2.] One defines the decreasing filtration 
$$
F^n \Omega^{\bullet}_{X, (E,e)} = \oplus_{a \geq n}
\Omega^{a,b}_{X, (E,e)}. 
$$
Then
$$
(\Omega^{\bullet}_{X}, \Omega^{\geq n}_{X}) \to
(\Omega^{\bullet}_{X, (E,e)}, F^n)
$$
is a filtered quasi-isomorphism. 
\item[3.] There is a canonical map 
$$
w^n: S^n (\sG^{*}_{e} )^{G_e} \to S^n (\sG^{*}_{(E,e)} ) =
\Omega^{n,n}_{X, (E,e)} 
$$
with ${\rm Im} \ w^n \subset (\Omega^{n,n}_{X, (E,e)})_{cl}
\subset F^n \Omega^{\bullet}_{X, (E,e)} [2 n]$ which defines
a map of commutative differential graded algebras 
$$
w: \oplus_n S^n (\sG^*_e )^{G_e} [-2n] \to \Omega^{\bullet}_{X,(E,e)} .
$$
This is called the {\it Weil homomorphism} of $(E,e)$ and
$\Omega^{\bullet}_{X,(E,e)}$ is called its {\it Weil algebra}.
(Locally on $X$, $T_{\sE_e /X} \simeq \sG_e \times G_e \times X$.
The gluing functions act via the adjoint representation on $\sG_e$
and the multiplication on $G_e$. Therefore $(S^n
\sG^{*}_{e})^{G_e}$ defines functions on $T_{\sE_e/X}$. This
defines $w^n$. To see that $d' ({\rm Im} \ w^n) = d'' ({\rm Im}
\ w^n) =0$, one may assume that $\sE_e \simeq G_e \times X$, and
since ${\rm Im} \ w^n$ does not depend on $X$, one also may
assume that $\sE_e \simeq G_e$. Then one knows that $d' : S^n
\sG^{*}_{e} \to \sG^{*}_{e} \otimes S^n \sG^{*}_{e}$ vanishes on
$(S^n \sG^{*}_{e})^{G_e}$, and $d'' |S^n \sG^{*}_{e} =0$ by
definition). In particular, $w^n$ defines a map 
$$
w^n : S^n (\sG^{*}_{e} )^{G_e} \to \H^{2n} (X, F^n
\Omega^{\bullet}_{X, (E,e)}) = \H^{2n} (X, \Omega^{\geq n}_{X}).
$$
\item[4.] When $k = \C$, one sets 
$$
S^n (\sG^{*}_{e})^{G_e}_{\Z (n)} = \Ker S^n (\sG^{*}_{e})^{G_e} \to
H^{2n} (X_{\an}, \C/\Z (n)). 
$$
For $E_e = E_{un,e}$, the universal bundle $E_{un,e} = G^{\Delta_
\ell}_{e} \times k^r /G_e$ on the simplicial $BG_e = G^{\Delta_
\ell}_{e} / G_e$, $w^n$ is an isomorphism, and one has 
$$
w^n : S^n ( \sG^{*}_{e} )^{G_e} \> \sim >> \H^{2n} (BG_e,
\Omega^{\geq n}) = \H^{2n} (BG_e, \Omega^{\bullet}), 
$$
the last isomorphism coming from 
$$
H^i (BG_e, \Omega^j ) = 0 \ i \neq j. 
$$
This remains true on $BG_{e,\an}$. This implies that $S^n
(\sG^{*}_{e})^{G_e}_{\Z(n)} \> \sim >> H^{2n} (BG_{e, \an} , \Z
(n))$, the group $H^{2n} (BG_{e,\an}, \Z (n))$ being torsion
free. 
\item[5.] When $k = \C$, they define the {\it Weil cohomology}
by 
$$
U_{(E,e)} (n) = {\rm cone} \ (\Z (n) \oplus S^n
(\sG^{*}_{e})^{G_e} [-2n] \to \Omega^{\bullet}_{X, (E,e)}) [-1],
$$
with the exact sequence 
\begin{multline*}
0 \to H^{2n-1} (X_{\an}, \C / \Z (n)) \to \\
H^{2n} (X_{\an},
U_{(E,e)} (n)) \to S^n (\sG^{*}_{e})^{G_e}_{\Z(n)} \to 0 .
\end{multline*}
The complexes $U_{(E,e)} (n)$ have a multiplication defined as
in \cite{B}, \cite{EV} for the Deligne cohomology, compatible
with the map of complexes
$$
U_{(E,e)} (n) \to {\rm cone} \ (Z(n) \oplus F^n
\Omega^{\bullet}_{X, (E,e)} \to \Omega^{\bullet}_{X, (E,e)} )
[-1] $$
induced by $w^n$. This induces a map, still denoted by $w^n$,
from the Weil cohomology to the analytic Deligne cohomology. 
\item[6.] A $\Omega^{1}_{X}$ valued connection $\nabla : (E,e)
\to \Omega^{1}_{X} \otimes (E,e)$ on $(E,e)$ is a splitting of
the $\Omega^{1}_{X, (E,e)}$ sequence, or, equivalently, a
quotient map of commutative differential graded algebras
$$
\nabla : \Omega^{\bullet}_{X, (E,e)} \to \Omega^{\bullet}_{X}
$$
which is a quasi-isomorphism inversing $\Omega^{\bullet}_{X} \to
\Omega^{\bullet}_{X, (E,e)}$. (Then $\nabla^2 =0$ if and only if
$\nabla$ induces a quasi-isomorphism of filtered complexes).
Thus $\nabla$ maps $U_{(E,e)} (n)$ to 
$$
C (n)^{\an} |X = {\rm cone} \ (\Z (n) \oplus \Omega^{\geq
2n}_{X} \to \Omega^{\bullet}_{X} ) [-1] 
$$
via 
$$
S^n (\sG^{*}_{e})^{G_e} [-2n] \> w^n >> F^{2n}
\Omega^{\bullet}_{X, (E,e)} \> \nabla >> \Omega^{\geq 2n}_{X}.
$$
This defines classes of $(E, \nabla,e)$ in 
$$
\H^{2n} (X_{\an}, C (n)^{\an} |X)
$$
(\ref{prop:diag}, 2)), mapping
to the classes in the analytic Deligne cohomology
(but not in $D^n(X)$ if $X$  
is not proper). The image of this class is $P_n (\nabla^2
, \ldots , \nabla^2)$ (see \ref{sec:not}, 10)), since $\nabla$
maps $\sG^{*}_{(E,e)}$ to $\Omega^{2}_{X}$ via the curvature. In
particular, if $\nabla^2 =0$, then those classes are in
$H^{2n-1} (X_{\an} , \C / \Z (n))$. The authors claim that this
is exactly $\hat{c}_n (E, \nabla)$ (without saying why) for $e
=0$. 
\end{enumerate}

\subsection{Algebraic construction}

We want to make the construction of Beilinson-Kazhdan algebraic.

\begin{rem} \label{rmk}
\begin{enumerate}
\item[1)] The construction of the filtered quasi-isomorphism 
$$
(\Omega^{\bullet}_{X} , \Omega^{\geq n}_{X}) \to
(\Omega^{\bullet}_{X,E}, F^n \Omega^{\bullet}_{X,E})
$$
is algebraic, as well as the Weil morphism $w^n$. If $(E,
\nabla, \tau_0)$ is as in \ref{sec:not}, 4), 6), one first
defines $\Omega^{\bullet}_{X,E} (\log D)$ by replacing the
Atiyah extension by the logarithmic Atiyah extension 
$$
0 \to \Omega^{1}_{X} (\log D) \to \Omega^{1}_{X,E} (\log D) \to
\sG^{*}_{E} \to 0.
$$
Then one defines $\Omega^{\bullet}_{X,E,\tau_0} (\log D)$ by
replacing the logarithmic Atiyah extension by the $N^1$ valued
Atiyah extension
$$
0 \to N^1 \to \Omega^{1}_{X, E, \tau_0} (\log D) \to \sG^{*}_{E}
\to 0
$$
obtained by push forward through $\Omega^{1}_{X} (\log D) \to
N^1$ of the logarithmic extension. Thus 
$$
\Omega^{n}_{X, E, \tau_0} (\log D) = \frac{N^n \oplus
\Omega^{n}_{X,E} (\log D)}{\Omega^{n}_{X} (\log D)}
$$
This defines a filtered quasi-isomorphism 
$$
(N^{\bullet}, N^{\geq n} ) \to (\Omega^{\bullet}_{X, E, \tau_0}
(\log D), F^n \Omega^{\bullet}_{X, E, \tau_0} (\log D)) 
$$
where $F^n$ is defined as in \ref{subsection:BK} 2): 
\begin{gather}
F^n \Omega^{\bullet}_{X,E, \tau_0} (\log D) = \oplus_{a \geq n}
\Omega^{a,b}_{X,E,\tau_0} (\log D) \notag \\
\Omega^{a,b}_{X,E,\tau_0} (\log D) = \Lambda^{a-b}
\Omega^{1}_{X,E,\tau_0} (\log D) \otimes S^b (\sG^{*}_{E} ).
\notag 
\end{gather}
We still denote by $w^n$ the induced Weil homomorphism
\begin{multline*}
(S^m \sG^{*})^G \to S^n (\sG^{*}_{E}) \to
(\Omega^{n,n}_{X,E,\tau_0} (\log D))_{cl} \\
\to F^n \Omega^{1}_{X,E,\tau_0} (\log D) [2n].
\end{multline*}
When $\pi : Y \to X$ is a morphism of (simplicial) smooth
varieties, such that $D' = \pi^{-1} (D)$ is a normal crossing
divisor, then one has a map 
$$
\pi^{-1} : \Omega^{\bullet}_{X,E,\tau_0} (\log D) \to R \pi_*
\Omega^{\bullet}_{Y, \pi^* E, \tau_0} (\log D')
$$
induced by
$$
\pi^{-1} : \Omega^{\bullet}_{X,\tau_0} (\log D) \to R \pi_*
\Omega^{\bullet}_{Y, \tau_0} (\log D')
$$
(\ref{dfn:tauform}). 
\item[2)] Then a $ N^1$ valued connection on $E$ defines a
quotient of commutative graded algebras
$$
\nabla: \Omega^{\bullet}_{X,E,\tau_0} (\log D) \to N^{\bullet}
$$
which is an inverse to the quasi-isomorphism 
$$
N^{\bullet} \to \Omega^{\bullet}_{X, E, \tau_0} .
$$
Then $\nabla^2 =0$ is equivalent to $\nabla$ being a filtered
quasi-isomorphism. 
\end{enumerate}
\end{rem}

\begin{defn} \label{dfn:ABK}

Let $E$ be a rank $r$ vector bundle and $\tau_0$ be as in
\ref{sec:not} 4). We define 
\begin{multline*}
ABK^{n}_{E, \tau_0} = {\rm cone} \ (\sK_n \otimes S^n
(\sG^{*})^G [-n] \\
\> d \log \oplus - w^n >> F^n \Omega^{\bullet}_{X,E,\tau_0}
(\log D) [n]) [-1] .
\end{multline*}
For $\tau_0 =$ identity, $D \neq \phi$ we write simply
$ABK^{n}_{E}$. 
\end{defn} 

\begin{lem} \label{lemma:un}

For $X = BG,E = E_{un}$ (\ref{subsection:BK} 4)) one has 
$$
\H^n (BG, ABK^{n}_{E_{un}} ) = H^n (BG, \sK_n) .
$$
For $X = BG_e$, and $E_{un,e} = E_{un}|_{BG_e}$ with
its canonical filtration coming from this flag, one has 
$$
\H^n (BG_e, ABK^{n}_{E_{un,e}} ) = H^n (BG_e, \sK_n) .
$$
\end{lem} 

\begin{proof} 
We just write the exact sequence
\begin{multline*}
\to \H^{2n-1} (X, N^{\geq n}) \to \H^n (X, ABK^{n}_{E, \tau_0,
e}) \\
\to H^n (X, \sK_n) \oplus S^n (\sG^{*}_{e})^G \to \H^{2n} (X,
N^{\geq n}) 
\end{multline*}
and apply \ref{subsection:BK}, 4). 
\end{proof}

\begin{defn}

Let $(E, \tau_0)$ be as in \ref{dfn:ABK}. We define 
$$
c^{ABK}_{n} (E) \in \H^n (X, ABK^{n}_{E, \tau_0} (\log D))
$$
as the inverse image under a map $[E]: X_{\bullet} \to BG$
defined on a simplicial model of $X$ by the transition functions
of $E$, of 
$$
c_n (E_{un}) \in H^n (BG, \sK_n) = \H^n (BG, ABK^{n}_{E_{un}}). 
$$
\end{defn}

As usual, all possible classes of $[E]$ are homotop and this
definition does not depend on the choice of $[E]$. However, one
has to apply here the functoriality of the complexes
$ABK^{n}_{E, \tau_0} (\log D)$ as explained in \ref{rmk} 1).

Let $(E, \nabla ,\tau_0)$ be now as in \ref{sec:not}, 4), 6).
Then $\nabla$ induces a map, still denoted by $\nabla$, from
$ABK^{n}_{E, \tau_0} (\log D)$ to 
$$
C(n)_{\tau_0} = {\rm cone} \ (\sK_n \oplus N^{\geq 2n} [n] \to
N^{\geq n} [n]) [-1]
$$
(\ref{rmk} 2)). 

\begin{defn} 
$$
c^{ABK}_{n} (E, \nabla) = \nabla c^{ABK}_{n} (E) 
\in AD^{n}_{\tau_0} (X). 
$$
\end{defn}

\begin{thm} \label{thm:comp}

One has 
$$
c^{ABK}_{n} (E, \nabla) = c_n (E, \nabla) \in AD^{n}_{\tau_0}
(X) .
$$
\end{thm}

\begin{proof}

We want to apply the characterization of the classes
\ref{thm:char}. 
\begin{enumerate}
\item[1)] The functoriality of the class $c^{ABK}_{n} (E)$ is
clear. This implies the functoriality of the classes
$c^{ABK}_{n} (E, \nabla)$. 
\item[2)] Let $E$ be of rank 1. Then $\sG_{E} = \sO_X$, and
$(\sG^*)^G = k$. One has 
$$
\H^1 (X, ABK^{1}_{E}) = \H^1 (X, \sK_1 \> d \log >>
\Omega^{1}_{X,E} \>>> \Lambda^{2} \Omega^{1}_{X,E} \oplus
\frac{\sO_X}{k} ). 
$$
One considers the principal $\G_m$ bundle $p: \sE = E - \{ 0 \}
\to X$ with local trivialization $\sE|_{U_i} = \G_m \times U_i$,
and gluing $t_i = \xi_{ij} t_j$, $t_i$ parameter of $\G_m$ and
$\xi_{ij} \in \sO^{*}_{X} (U_i \cap U_j)$. This defines a
$\Omega^{1}_{X,E} = (p_* \Omega^{1}_{\sE})^{\G_m}$ connection on
$E$ via $\frac{d t_i}{t_i} - \frac{d t_j}{t_j} = \frac{d
\xi_{ij}}{\xi_{ij}}$, with local form $\frac{d t_i}{t_i}$. Then
$d'(\frac{dt_i}{t_i}) = 0$, while $d'' (\frac{dt_i}{t_i})$, which
is the image of $\frac{dt_i}{t_i}$ in $\sO_X$, is the residue of
$\frac{dt_i}{t_i}$ along the zero section of $E$. Therefore,
this is $1 \in k \subset \sO$. Thus $(d' + d'')
(\frac{dt_i}{t_i}) = 0$ in $\frac{\sO_X}{k}$ and $(\xi_{ij},
\frac{dt_i}{t_i})$ defines the class of $E$ in $\H^1 (X,
ABK^{1}_{E})$. Now if $\nabla$ is a $\tau_0$ connection on $E$,
it maps $\frac{dt_i}{t_i}$ to the local form $\alpha_i$ of the
connection, Thus 
$$
\nabla (\xi_{ij}, \frac{dt_i}{t_i}) = (\xi_{ij}, \alpha_i ) =
c_1 (E, \nabla) \in AD^{1}_{\tau_0} (X).
$$
\item[3)] In order to understand the Whitney product formula,
one has first to introduce a product on $ABK^{n}_{E, \tau_0}
(\log D)$. Again, one takes the same definition as in
\ref{dfn:pro}: 
$$
ABK^{m}_{E,\tau_0} (\log D) \times ABK^{n}_{E,\tau_0} (\log D)
\> \cup_{\alpha} >> ABK^{m+n}_{E, \tau_0} (\log D)
$$
is given by
$$
\begin{array}{lll}
x \cup_{\alpha} y & = \{ x,y\} & x \in \sK_m , y \in \sK_n \\
& & \\
& = 0 & x \in \sK_m, y \in S^n \\
& & \\
& = (1- \alpha) d \log x \wedge y & x \in \sK_m, y \in F^n \\
& & \\
& = 0 & x \in S^m, y \in \sK_n \\
& & \\
& = x \wedge y & x \in S^m, y \in S^n \\
& & \\
& = (-1)^{{\rm deg} \ x} \alpha x \wedge y & x \in S^m, y \in
F^n \\
& = \alpha x \wedge d \log y & x \in S^m, y \in \sK_n \\
& & \\
& = (1-\alpha) x \wedge y & x \in F^m , y \in S^n \\
& & \\
& = 0 & x \in F^m , y \in F^n , 
\end{array} $$
where we shortened the notations by 
$$
\begin{array}{ll}
S^m = & S^m (\sG^*)^G [-m] \\
F^m = & F^m \Omega^{\bullet}_{X,E,\tau_0} (\log D) [m].
\end{array} 
$$
Then a $N^1$ valued connection $\nabla$ maps $ABK^{n}_{E,
\tau_0} (\log D)$ to $C(n)_{\tau_0}$ (\ref{rmk} 2)) compatibly
to the product.
\item[4)] Let $e: 0 \to E'' \to E \to E' \to 0$ be an exact
sequence of bundles, with rank $E'' = r''$, rank $E' = r'$, $r =
r' +r''$. Let $\sG_{E,e} \>\rho>> \sG_E$ be the endomorphisms
of $E$ respecting the extension $e$, and $\sG_{E,e} \> p'
\oplus p'' >> \sG_{E'} \oplus \sG_{E''}$ be the restrictions of the
compatible endomorphisms to those of $E'$ and of $E''$. This defines maps
$$
\sG^{*}_{E} \> \rho^* >> \sG^*_{E,e} \< (p' \oplus
p^n)^* << \sG^{*}_{E'} \oplus \sG^{*}_{E''}
$$
and maps 
$$
ABK^{n}_{E} \> \rho^* >> ABK^{n}_{E,e} \< p^{'*} \oplus p^{''*}
<< ABK^{n}_{E'} \oplus ABK^{n}_{E''}
$$
where one defines $ABK^{n}_{E,e}$ in the following way. Let $G_e
\subset G$ be the matrices of the form $$\left( \begin{array}{ll}
A & C \\ 0 & B \end{array} \right), A \in GL( r''), B \in GKL
(r'),$$ and $\sG_e \subset \sG$ be the corresponding algebra.
Then by restriction of the structure group of $E$ from $G$ to
$G_e$, one considers $\sE_e$ the corresponding $G_e$ torseur and
one does the same construction as in \ref{subsection:BK} with
$\sE_e$ replacing $\sE$. Now $E_{un} |_{BG_e}$ is the extension
$$
e_{un} : 0 \to E''_{un} \to E_{un} |_{BG_e} \to E'_{un} \to 0
$$
and one has 
\begin{multline*}
H^n (BG, \sK_n) = \H^n(BG, ABK^n_{E_{un}}) \\
\> \rho^* >> H^n
(BG_e, ABK^n_{(E_{un}, e_{un})}) = H^n (BG_e , \sK_n) 
\end{multline*}
( Lemma \ref{lemma:un}) and 
$$
H^n (BG_e, \sK_n) = \oplus^{n}_{a=0} p^{'*} H^a (BG', \sK_a)
\cup p^{''*} H^{n-a} (BG'', \sK_{n-a}).
$$
Since the product on $ABK$ is compatible with the product on
$\sK$, this implies 
\begin{multline*}
H^n (BG_e, ABK^n_{(E_{un}, e_{un})}) = \\
\oplus^{n}_{a=0} p' \H^a
(BG', ABK^a_{E'_{un}}) \cup p'' \H^{n-a} (BG'', ABK^{n-a}_{E''_{un}})
\end{multline*}
and one has the decomposition
$$
\rho^* c^{ABK}_{n} (E_{un}) = \oplus^{n}_{a=0} p^{'*}
c^{ABK}_{a} (E'_{un}) \cup p^{''*} c^{ABK}_{n-a} (E''_{un}). 
$$
By functoriality, one obtains the similar relation for the
classes of $E$: 
$$
\rho^* c^{ABK}_{n} (E) = \oplus^{n}_{a=0} p^{'*} c^{ABK}_{a}
(E') \cup p^{''*} c^{ABK}_{n-a} (E'').
$$
If $E$ now has a $\tau_0$ connection $\nabla$, the map $\nabla :
\Omega^{\bullet}_{X,E,\tau_0} (\log D) \to N^{\bullet}$
(\ref{rmk} 2)) factors through
$$
\Omega^{\bullet}_{X,E, \tau_0} (\log D) \>\rho^* >>
\Omega^{\bullet}_{X,E,e,\tau_0} (\log D) \> \nabla_{\rho} >>
N^{\bullet} 
$$
when $\nabla$ is compatible with the exact sequence. This shows
the Whitney product formula. 
\end{enumerate}
\end{proof}

\begin{rem}
The construction of the tautological $\Omega^{1}_{X,E}$ valued
connection on $E$ of rank $r$ goes as in the second section of
the proof of \ref{thm:comp} for $r =1$: One considers a local
trivialization $\sE |_{U_i} \simeq G \times U_i$ of the
principal $G$ bundle $\sE$. Then the gluing is given by $g_i =
g_{ij}g_{j}$, $g_{ij} \in G (U_i \cap U_j)$, where $g_i$ is
the tautological $G$ valued function on $G : g_i (x) = x$. Then 
$$
dg_i g^{-1}_{i} - g_{ij} dg_j g_j^{-1}g^{-1}_{ij} = dg_{ij} g^{-1}_{ij}
$$
is the equation of the connection. 
\end{rem}\bibliographystyle{plain}
\renewcommand\refname{References}
 
\end{document}